\documentclass[12pt]{article}
\usepackage{amssymb, amsmath}
\begin{document}
\renewcommand{\theequation}{\arabic{section}.\arabic{equation}}
\newtheorem{theorem}{Theorem}[section]
\newtheorem{lemma}{Lemma}[section]
\newtheorem{definition}{Definition}
\newtheorem{pro}{Proposition}[section]
\newtheorem{cor}{Corollary}[section]
\newcommand{\n}{\nonumber}
\newcommand{\tv}{\tilde{v}}
\newcommand{\tw}{\tilde{\omega}}
\renewcommand{\t}{\theta}
\newcommand{\w}{\omega}
\newcommand{\e}{\varepsilon}
\renewcommand{\a}{\alpha}
\renewcommand{\l}{\lambda}
\newcommand{\vare}{\varepsilon}
\newcommand{\s}{\sigma}
\renewcommand{\o}{\omega}
\renewcommand{\O}{\Omega}
\newcommand{\bb}{\begin{equation}}
\newcommand{\ee}{\end{equation}}
\newcommand{\bq}{\begin{eqnarray}}
\newcommand{\eq}{\end{eqnarray}}
\newcommand{\bqn}{\begin{eqnarray*}}
\newcommand{\eqn}{\end{eqnarray*}}
%  Cause equations to be numbered within sections.
\title{Nonexistence of  self-similar singularities in
the viscous magnetohydrodynamics  with zero resistivity}
\author{Dongho Chae\thanks{This work was supported partially by  KRF Grant(MOEHRD, Basic
Research Promotion Fund) and the KOSEF Grant no.
R01-2005-000-10077-0.}\\
Department of Mathematics\\
              Sungkyunkwan University\\
               Suwon 440-746, Korea\\
              {\it e-mail : chae@skku.edu}}
 \date{}
\maketitle
\begin{abstract}
We are concerned on the possibility of finite time singularity in a
partially viscous magnetohydrodynamic equations in $\Bbb R^n$,
$n=2,3$, namely the MHD with positive viscosity and zero
resistivity. In the special case of zero magnetic field the system
reduces to the  Navier-Stokes equations in $\Bbb R^n$.
 In this paper we exclude the scenario of finite time singularity in
the form of self-similarity, under suitable integrability conditions
on the velocity and the magnetic field. We also prove the
nonexistence of asymptotically self-similar singularity. This
provides us information on the behavior of solutions near possible
singularity of general type as described in Corollary 1.1 below.
\end{abstract}
\section{Introduction}
 \setcounter{equation}{0}
The equations of magnetohydrodynamics(MHD) with zero resistivity in
$\Bbb R^n$, $n=2,3$, are the following.
  \bq \label{1.1}
 &&\frac{\partial v}{\partial t} +(v\cdot \nabla )v =
 \nu \Delta v-\nabla (p +\frac12 |b|^2)+(b\cdot\nabla)b , \\
 \label{1.2}
 &&\frac{\partial b}{\partial t} +(v\cdot \nabla )b =(b \cdot \nabla )
 v,\\
 \label{1.3}
 &&\textrm{div }\, v =\textrm{div }\, b= 0 ,\\
 \label{1.4}
  &&v(x,0)=v_0 (x), \quad b(x,0)=b_0 (x),
  \eq
where $v=(v_1, \cdots , v_n )$, $v_j =v_j (x, t)$, $j=1,\cdots,n$,
is the velocity of the flow, $b=(b_1, \cdots , b_n )$, $b_j =b_j (x,
t)$, is the magnetic field, $p=p(x,t)$ is the scalar pressure,
$\nu>0$ is the viscosity of the fluid, and $v_0$, $b_0$ are the
given initial velocity and magnetic fields,
 satisfying div $v_0 =\mathrm{div}\, b_0= 0$, respectively.
 The system (\ref{1.1})-(\ref{1.4}) describes the macroscopic behavior of electrically
 conducting incompressible fluids with extremely high conductivity.
In the original (fully viscous) equations of magnetohydrodynamics,
besides the viscosity term, $\nu \Delta v$, in (\ref{1.1}) we have
the resistivity term, $\eta\Delta b$, in the right hand side of
(\ref{1.2}), where $\eta$ is
 the resistivity constant, which  is inversely proportional to
  the electrical conductivity constant, $\sigma$.
 In the extremely high electrical conductivity cases, which occur
 frequently in the cosmical and geophysical problems we ignore the
 resistivity term to have our system (\ref{1.1})-(\ref{1.4})(see e.g. \cite{chan}). We are
 concerned here the mathematical question of the global
 well-posedness/finite time singularity of the system
 (\ref{1.1})-(\ref{1.4}). The proof of local well-posedness of the Cauchy problem is rather
standard(actually the necessary essential estimates are derived in
the proof of Lemma 2.1 below), and similar to the case of fully
viscous MHD is done in \cite{tem}. The question of spontaneous
apparition of
 singularity  from a local classical solution is a challenging open problem
 in the mathematical fluid mechanics. The situation is similar to
 the  both of the cases of ideal MHD and
 fully viscous MHD.
 We just refer some of the studies on the finite time blow-up problem in the
 ideal MHD(\cite{caf,cor1,cor2, gra, gib} and references therein).
 In order to discuss the self-similar singularity of the system
 (\ref{1.1})-(\ref{1.4})
 we first observe that  it has the
following  scaling property: If $(v, b, p)$ is a solution of
(\ref{1.1})-(\ref{1.4}) corresponding to the initial data $(v_0,
b_0)$, then for any $\lambda
>0$ the functions
$$
  v^{\lambda}(x,t)=\lambda v (\lambda x, \l^{2}
  t), \quad  b^{\lambda}(x,t)=\lambda b (\lambda x, \l^{2}
  t),$$
 and
  $$ \quad p^{\l}(x,t)=\l^{2}p(\l x, \l^{2} t )
$$
  are also solutions with the initial data
  $ v^{\lambda}_0(x)=\lambda  v_0
   (\lambda x)$,  $ b^{\lambda}_0(x)=\lambda  b_0
   (\lambda x)$.
 In view of the above scaling
  property the  self-similar blowing up
  solution $(v(x,t), b(x,t))$ of the system (\ref{1.1})-(\ref{1.4}), if it exists,
   should be of the form,
  \bq
  \label{vel}
 v(x, t)&=&\frac{1}{\sqrt{T_*-t}}
V\left(\frac{x}{\sqrt{T_*-t}}\right),\\
  \label{mag}
 b(x, t)&=&\frac{1}{\sqrt{T_*-t}}
B\left(\frac{x}{\sqrt{T_*-t}}\right),\\
\label{pre}
  p(x,t)&=&\frac{1}{T_*-t}
P\left(\frac{x}{\sqrt{T_*-t}}\right)
 \eq
 for $t$
 close to the possible blow-up time $T_*$.
 If we substitute (\ref{vel})-(\ref{pre})
 into (\ref{1.1})-(\ref{1.4}), then we find that $(V,B,P)$ should be a solution of the
 stationary system:
\bq
 \label{le1} && \frac{1}{2}V +\frac{1}{2} (y\cdot \nabla)
V+(V\cdot \nabla)
 V=\nu \Delta V+(B\cdot \nabla )B \n\\
 &&\qquad \hspace{2.5in}-\nabla (P +\frac12 |B|^2), \\
 \label{le2}
&&\frac{1}{2}B+\frac{1}{2} (y\cdot \nabla)B+(V\cdot
 \nabla )B=(B\cdot \nabla )V,\\
 \label{le3}
&&\mathrm{ div}\, V=\mathrm{div}\, B=0.
 \eq
 Conversely, if $(V, B, P)$ is a smooth solution of the system (\ref{le1})-(\ref{le3}),
 then the triple of functions $(v,b,p)$ defined  by (\ref{vel})-(\ref{pre}) is
 a smooth solution of (\ref{1.1})-(\ref{1.4}),
  which blows up at $t=T_*$.
 The search for self-similar singularities of the form, (\ref{vel})-(\ref{pre})
was suggested first by Leray  for
 the 3D Navier-Stokes
 equations in \cite{ler}, and its nonexistence was first proved by Ne$\check{c}$as,
 Ru$\check{z}$i$\check{c}$ka
  and $\check{S}$ver$\acute{a}$k in \cite{nec} under the condition of
   $V\in L^3(\Bbb R^3)\cap H^1_{loc}
 (\Bbb R^3 )$, the result of which was generalized later by Tsai
 to the case $L^p (\Bbb R^3 )\cap H^1 _{loc} (\Bbb R^3)$ with $p>3$ in
 \cite{tsa}. Their proofs crucially depend on the maximum principle
 of the Leray system,
$$
\frac12 V+\frac12 (y\cdot \nabla) V+(V\cdot \nabla)
 V=-\nabla P +\nu \Delta V, \quad \mathrm{div}\,V=0,
 $$
which corresponds to a special case($B=0$) in
(\ref{le1})-(\ref{le3}). The corresponding maximum principle for
(\ref{le1})-(\ref{le3}), however, cannot be obtained by applying
similar method used in \cite{nec, tsa}(The situation is similar even
if we have `special' resistivity term $\nu \Delta B$ to the right
hand side of (\ref{le2})).
  Due to this fact there are difficulties in extending the
 nonexistence results for the self-similar singularity
 of the 3D Navier-Stokes system to our system (\ref{1.1})-(\ref{1.4}).
  Recently, the author of this paper developed
 new method to prove nonexistence of the self-similar singularity of
 the 3D Euler system under suitable integrability condition on the
 vorticity(\cite{cha1}).
Here we first combine the argument in (\cite{cha1}) together with
the  results by \cite{nec, tsa} to obtain the nonexistence of
self-similar blowing up solutions, the precise statement of which is
in the
 following theorem.
\begin{theorem}
Suppose there exists $T_*>0$ such that we have  a
 representation of a solution $(v, b)$
  to  (\ref{1.1})-(\ref{1.4})by
  (\ref{vel})-(\ref{mag})  for all $t \in (0 ,T_*)$   with $(V, B)$
  satisfying the following
  conditions:
\begin{itemize}
 \item[(i)] $(V, B)\in [C ^1 (\Bbb R^n)]^2$,  $\nabla V\in
 L^\infty(\Bbb R^n )$, and  div $V=$ div $B=0$.
 \item[(ii)] In the case $n=3$, there exists $q_1\in [3, \infty)$
 such that $V\in L^{q_1}(\Bbb R^3)$. \\
 In the case $n=2$, $V\in L^2 (\Bbb R^2)$.
  \item[(iii)] There exists $q_2 >0$ such that
 $B\in L^q (\Bbb R^n
  )$ for all $q\in (0, q_2)$.
\end{itemize}
  Then, $V=B=0$.
\end{theorem}
\noindent{\textsf{Remark} 1.1} In order to illustrate the
integrability condition for $B$ in (iii) above  we make the
following observations: If a function $f(x)$ on $\Bbb R^n$ satisfies
$$
\sup_{x\in \Bbb R^n}(1+|x|^k) |f(x)|<\infty \quad \forall k\in \Bbb
N,
$$
then $f\in L^p (\Bbb R^n)$ for all $p\in (0, \infty)$. Indeed, given
$p\in (0, \infty)$, we choose $k=[\frac{n+1}{p}]$. Then, we have
$$ \int_{\Bbb R^n} |f(x)|^p dx
\leq  \int_{\Bbb R^n} \left( \frac{C}{1+|x|^k}\right)^{p} dx\leq
 C(n,p) \int_0 ^\infty \frac{r^{n-1}}{(1+r^{n+1}) }dr <\infty.
 $$
\ \\
Under different type of decay conditions on $(V, B)$ from the above
theorem, we could also have similar nonexistence result as follows.
\begin{theorem}
Suppose there exists $T_*>0$ such that we have  a
 representation of a solution $(v, b)$
  to (\ref{1.1})-(\ref{1.4}) by
  (\ref{vel})-(\ref{mag})  for all $t \in (0 ,T_*)$   with $(V, B)$
  satisfying the following
  conditions:\\
\begin{itemize}
 \item[(i)]$(V, B) \in [H^m (\Bbb R^n)]^2$, $m>n/2+1$
\item[(ii)] $\|\nabla {V}\|_{L^\infty} +\|\nabla
{B}\|_{L^\infty}<\eta$, where $\eta$ is a sufficiently small
constant to be determined in Lemma 2.1 in the next section.
\end{itemize}
  Then,  ${V}={B}=0$.
\end{theorem}
\noindent{\textsf{Remark 1.2}} The above theorem  implies the {\it
`stability of the null solution'} of  the stationary system
(\ref{le1})-(\ref{le3}). Namely, there exists $\eta>0$ such that if
$(V, B)$ is a solution to (\ref{le1})-(\ref{le3}) and belongs to a
ball $B(0,\eta)=\{ X=(V, B)\in H^m (\Bbb R^n ) \, ;\, \|\nabla
X\|_{L^\infty} <\eta \}$, where $m>n/2+1$, then $(V,B)=(0,0)$. \\
\ \\
Next, we consider  more refined scenario of `asymptotically
self-similar singularity', which
 means that the local in time smooth solution evolves into a
self-similar profile as the possible singularity time is approached.
 A similar notion was
 considered previously by Giga and Kohn in the context of the
 nonlinear scalar heat equation in \cite{gig}.
 Recently, the author of this paper
 (\cite{cha2}) considered it in the context of 3D Euler and the 3D
 Navier-Stokes equations(see also
 \cite{hou}), and excluded its scenario.
We apply the idea developed in \cite{cha2} to exclude asymptotically
self-similar singularity of our system (\ref{1.1})-(\ref{1.4}).
 \begin{theorem}
Let $(v, b) \in [C([0, T);H^m (\Bbb R^n))]^2$, $m>n/2+1$,  be a
classical solutions to (\ref{1.1})-(\ref{1.4}).
 Suppose there exist  functions $\bar{V}, \bar{B}$ satisfying the conditions (i)-(iii) for
 $V, B$ in Theorem 1.1 such that the following boundedness and the convergence hold true:
  \bq\label{th12}
  \lefteqn{\sup_{0<t<T}
(T-t)^{\frac{1-n}{2}}\left\|v(\cdot, t) -\frac{1}{\sqrt{T-t}}
\bar{V} \left(\frac{\cdot}{\sqrt{T-t}} \right)\right\|_{L^1
}}\hspace{.0in}\n \\
&&+ \sup_{0<t< T} (T-t)^{\frac{1-n}{2}}\left\|b(\cdot, t)
-\frac{1}{\sqrt{T-t}} \bar{B} \left(\frac{\cdot}{\sqrt{T-t}}
\right)\right\|_{L^1 }
<\infty,\n \\
  \eq
  \bq\label{th122}
   \lefteqn{\lim_{t\nearrow T}
 (T-t)  \left\|\nabla v(\cdot, t) -\frac{1}{ \sqrt{T-t} } \nabla \bar{V}
\left(\frac{\cdot}{\sqrt{T-t}}
\right)\right\|_{L^\infty }}\hspace{.0in}\n \\
 &&+
 \lim_{t\nearrow T}
(T-t)\left\|\nabla b(\cdot, t) -\frac{1}{T-t}  \nabla \bar{B}
\left(\frac{\cdot}{\sqrt{T-t}}
\right)\right\|_{L^\infty} =0.\n \\
  \eq
 Then, $\bar{V}=\bar{B}=0$, and $(v, b)$ can be extended to a solution of
(\ref{1.1})-(\ref{1.4}) in $[0, T+\delta]\times \Bbb R^n$,
 and belongs to $C([0, T+\delta]; H^m (\Bbb R^n))$ for some $\delta
 >0$.
\end{theorem}
\noindent{\textsf{Remark} 1.3} Unlike to the cases of the Euler
equations(\cite{cha2}), the convergence of (\ref{th122}) is not in
the critical Besov space norms for the quantities of vorticities and
current densities, but in the Lipshitz norm for the gradients of
velocities and  magnetic fields. Actually due to the non-symmetry of
the viscosity terms(the term $\nu \Delta v$ for the velocity
evolution equations (\ref{1.1}), and zero for the magnetic field
evolution equations (\ref{1.2})) we cannot obtain critical Besov
space type of norm estimates in the procedure of proof of the above
theorem(see
the proof in the next section below).\\
 \ \\
As an immediate corollary of Theorem 1.3 we have the following
information of the behaviors of solution near possible singularity,
which is not necessarily of the self-similar type.
\begin{cor}
Let $(v, b) \in [C([0, T_*);H^m (\Bbb R^n))]^2$, $m>n/2+1$,  be a
classical solutions to (\ref{1.1})-(\ref{1.4}), which blows up at
$T$.  We expand the solution of the form:
 \bq
 v(x,t)&=&\frac{1}{\sqrt{T-t}}
\bar{V} \left(\frac{x}{\sqrt{T-t}}\right) +\bar{v}(x,t),\\
b(x,t)&=&\frac{1}{\sqrt{T-t}} \bar{B}
\left(\frac{x}{\sqrt{T-t}}\right) +\bar{b}(x,t),
 \eq
 where $(\bar{V}, \bar{B})$ satisfies the conditions (i)-(iii) for
 $(V, B)$ in Theorem 1.1.
 Then, either
 \bb
 \lim\sup_{t\nearrow T} \left[(T-t)^{\frac{1-n}{2}}(\|\bar{v}(t)\|_{L^1}
 +\|\bar{b}(t)\|_{L^1})\right]=\infty,
 \ee
 or there exists $\e_0 >0$ such that
 \bb
 \lim\sup_{t\nearrow T}\left[(T-t)(\|\nabla
 \bar{v}(t)\|_{L^\infty}
 +\|\nabla \bar{b}(t)\|_{L^\infty})\right] >\e_0.
 \ee
\end{cor}

\section{Proof of the theorems}

{\bf Proof of Theorem 1.1} We assume classical solution $(v, b)$ of
the form (\ref{vel})-(\ref{mag}). We will show that this assumption
leads to $V=B=0$. By consistency with the initial condition, $ b_0
(x)=\frac{1}{\sqrt{T_*}}B(\frac{x}{\sqrt{T_*}})$, we can rewrite the
representation (\ref{mag}) in the form,
 \bb \label{thm12aa}
 b (x,t) =\left(1-\frac{t}{T_*}\right)^{-\frac12} b_0
 \left( \left(1-\frac{t}{T_*}\right)^{-\frac12}x\right)
 \qquad \forall t\in [0, T_*).
 \ee
   Let $a\mapsto X(a,t)$ be the
particle
 trajectory mapping, defined by the ordinary differential equations,
 $$
 \frac{\partial X(a,t)}{\partial t} =v(X(a,t),t) \quad;\quad X(a,0)=a.
 $$
We set $A(x,t):=X^{-1} (x ,t)$, which is called the back-to-label
map, satisfying
 \bb\label{inverse}
 A (X(a,t),t)=a, \quad X(A
 (x,t),t)=x.
 \ee
 We note that for our smoothness condition (i) decay condition on the velocity (ii)
 the existence of $A(\cdot, t)$ is guaranteed at lease for $t$ close
 to $T_*$(see \cite{con}), which is enough for our purpose in the proof.
Taking dot product (\ref{1.2}) by $b$, we obtain
 \bb\label{new0}
 \frac{\partial |b|}{\partial t} +(v\cdot \nabla )|b |=\a |b|,
 \ee
where $\a (x,t)$ is defined as
 $$
 \a (x,t)=\left\{ \aligned &\sum_{i,j=1}^n S_{ij}
 (x,t)\xi_i(x,t)\xi_j (x,t) &\quad \mbox{if} \quad b (x,t)\neq 0\\
 &\qquad 0 &\quad \mbox{if} \quad b (x,t)=0\endaligned \right.
 $$
 with
 $$ S_{ij} =\frac12 \left( \frac{\partial v_j}{\partial x_i}
 +\frac{\partial v_i}{\partial x_j}\right)\quad\mbox{and}\quad
\xi (x,t)=\frac{b(x,t)}{|b (x,t)|}.
$$
 In terms of the particle trajectory mapping we
can rewrite (\ref{new0}) as
 \bb\label{new1}
 \frac{\partial }{\partial t} |b
(X(a,t),t)|=\alpha (X(a,t),t) |b (X(a,t),t)|.
 \ee
Integrating (\ref{new1}) along the particle trajectories $\{
X(a,t)\}$, we have
 \bb\label{new2}
 |b (X(a,t),t)|=|b_0 (a)|\exp \left[ \int_0 ^t \a
(X(a,s),s) ds \right].
  \ee
Taking into account the simple estimates
$$ -\|\nabla v(\cdot ,t)\|_{L^\infty}\leq  \a
(x,t)  \quad \forall x\in \Bbb R^n,
 $$
we obtain from (\ref{new2}) that
$$
|b_0 (a)|\exp \left[- \int_0 ^t \|\nabla v
      (\cdot,s)\|_{L^\infty} ds \right]\leq  |b (X(a,t),t)|,
$$
 which, in terms of the back-to-label map, can be rewritten as
  \bb\label{new3}
|b_0 (A(x,t))|\exp \left[- \int_0 ^t \|\nabla v
      (\cdot,s)\|_{L^\infty} ds \right]\leq |b(x,t)|.
 \ee
 Combining this with the self-similar representation formula in (\ref{thm12aa}), we have
 \bb
  \label{new4}
      |b_0 (A(x,t))|\exp \left[- \int_0 ^t \|\nabla v
      (\cdot,s)\|_{L^\infty} ds \right]\leq \left(1-\frac{t}{T_*}\right)^{-\frac12}
      |b_0\left( \left(1-\frac{t}{T_*}\right)^{-\frac12}x\right)|
 \ee
 Given $q\in (0, q_2)$, computing $L^q(\Bbb R^n)$ norm of the each side of (\ref{new4}),
 we obtain
 \bb\label{new4a}
   \|b_0 \|_{L^q} \exp \left[- \int_0 ^t \|\nabla v
      (\cdot,s)\|_{L^\infty} ds \right]\leq
       \|b _0 \|_{L^q} \left(1-\frac{t}{T_*}\right)^{\frac{n}{2q}-\frac12}
\ee
 where we used the fact $\mathrm{det} (\nabla A (x,t))\equiv 1$.
 Now, suppose $B\neq 0$, which is equivalent to assuming that $b_0 \neq
 0$, then we divide (\ref{new4a}) by $\|b_0 \|_{L^q}$ to have
 \bb \label{new5}
  \exp \left[- \int_0 ^t \|\nabla v
      (\cdot,s)\|_{L^\infty} ds \right]\leq  \left(1-\frac{t}{T_*}\right)^{\frac{n}{2q}-\frac12}
 \ee
Passing $q\searrow 0$ in
 (\ref{new5}), we deduce that
 $$\int_0 ^{t}\|\nabla v
      (\cdot,s)\|_{L^\infty} ds =\infty \quad \forall t\in (0, T_*)
      $$
      This contradicts with the assumption that the flow is smooth on $(0, T_*)$, i.e $v\in
      C^1 ([0,T_*); C^1 (\Bbb R^n)\cap W^{1, \infty}(\Bbb R^n))$, which is implied by the
by the explicit representation formula (\ref{vel})-(\ref{mag}),
combined with the assumption (i). Hence we need to have
      $B=0$. Setting  $B=0$ in the system (\ref{1.1})-(\ref{1.4}), it reduces to the
     incompressible Navier-Stokes system in $\Bbb R^n$. When $n=3$ we  apply
     Ne$\check{c}$as-Ru$\check{z}$i$\check{c}$ka-$\check{S}$ver$\acute{a}$k's
     result in
     \cite{nec} for $q_1=3$ and Tsai's result
      in \cite{tsa} for $q_1\in (3, \infty)$ respectively. Then,
     we obtain  $V=0$.  In the case $n=2$ we recall that
      in the 2D Navier-Stokes equations for the
     initial data $v_0 (\cdot)=\frac{1}{\sqrt{T_*}}V (\frac{\cdot}{\sqrt{T_*}} )
     \in L^2(\Bbb R^2 )$ the solution $v$ belongs to
     $C^\infty ((0, \infty)\times \Bbb R^2 )$(see e.g. \cite{tem0}),
     and hence we need to have $V=0$. $\square$\\
\ \\
In order to prove Theorem 1.2 and Theorem 1.3 we  establish the
following continuation principle for local classical solution of
(\ref{1.1})-(\ref{1.4}).
\begin{lemma}
Let $(v, b) \in [C([0, T);H^m (\Bbb R^n))]^2$, $m>n/2+1$,  be a
classical solution to (\ref{1.1}(-(\ref{1.4}). There exists an
absolute constant $\eta
>0$ such that if
 \bb\label{th11}
 \sup_{0\leq t<T} (T-t)
 \left\{\|\nabla v (t)\|_{L^\infty}
 +\|\nabla b(t)\|_{L^\infty }\right\}<\eta
 ,
 \ee
 then the solution $(v(x,t), b(x,t))$  can be extended to be
 functions on $[0, T+\delta]\times \Bbb R^n$,
 and belongs to $C([0, T+\delta]; H^m (\Bbb R^n))$ for some $\delta
 >0$.
\end{lemma}
 \noindent{\bf Proof} Let $\a =(\a_1, \cdots, \a_n
)\in (\Bbb N \cup \{ 0\})^n$ be a standard multi-index with $|\a
|=\a_1 +\cdots +\a_n$. We take operation $D^\a =\partial_1 ^{\a_1}
\cdots
\partial_n ^{\a_n}$ on (\ref{1.1}), and  take $L^2(\Bbb R^n )$ inner product it with
$D^\a$, summing over $|\a |\leq m$ after integration by parts. Then,
we obtain
 \bq\label{lemm1}
  \lefteqn{\frac12 \frac{d}{dt} \|v\|_{H^m} ^2 +\nu\|\nabla v\|_{H^m}^2
  =\sum_{|\a|\leq m} (D^\a (v\cdot \nabla )v -(v\cdot \nabla )D^\a v
  , D^\a v)_{L^2}}\hspace{.0in}\n \\
  && +\sum_{|\a|\leq m} (D^\a (b\cdot \nabla )b -(b\cdot \nabla )D^\a
  b , D^\a  v)_{L^2} +\sum_{|\a|\leq m} (b\cdot \nabla )D^\a b
  , D^\a v)_{L^2},\n\\
  \eq
  where we used the facts,
  $$\left((v\cdot \nabla )D^\a v
  , D^\a v\right)_{L^2}=\frac12 \int_{\Bbb R^n} (v\cdot \nabla )|D^\a v|^2
  dx=-\frac12 \int_{\Bbb R^n } (\mathrm{div}\, v) |D^\a v|^2 dx =0,
  $$
  and
  $$(D^\a v, D^\a \nabla (p+\frac12 |v|^2 ) )_{L^2}=-(D^\a(\mathrm{div}\,
  v), D^\a(p+\frac12 |v|^2 )  )_{L^2}=0.
  $$
  Applying the well-known commutator estimate(\cite{kla}),
  $$
   \sum_{|\a|\leq m}\|D^\a (fg)-fD^\a g\|_{L^2} \leq C (\|\nabla
   f\|_{H^{m-1}} +\|f\|_{H^m} \|g\|_{L^\infty} ),
   $$
   to the terms of the right hand side of (\ref{lemm1}), we have
   \bq\label{lemm2}
   \lefteqn{\frac12 \frac{d}{dt} \|v\|_{H^m} ^2 +\nu\|\nabla v\|_{H^m}^2 \leq
   C \|\nabla v \|_{L^\infty} \|v\|_{H^m} ^2}\hspace{0.in}\n \\
   && +C \|\nabla b
   \|_{L^\infty} \|b\|_{H^m} \|v\|_{H^m} + \sum_{|\a|\leq m}((b\cdot \nabla )D^\a b, D^\a v
   )_{L^2}.
   \eq
   Similarly, starting from (\ref{1.2}), we can deduce
   \bq\label{lemm3}
  \lefteqn{\frac12 \frac{d}{dt} \|b\|_{H^m} ^2=
  -\sum_{|\a|\leq m} (D^\a (v\cdot \nabla )b -(v\cdot \nabla )D^\a b
  , D^\a b)_{L^2}}\hspace{.0in}\n \\
  &&+\sum_{|\a|\leq m} (D^\a (b\cdot \nabla )v -(v\cdot \nabla )D^\a
  v , D^\a b)_{L^2} +\sum_{|\a|\leq m}((b\cdot \nabla )D^\a v, D^\a
  b )_{L^2}\n \\
  &&\leq C\|\nabla v\|_{L^\infty} \|b\|_{H^m}^2 +C \|\nabla
  b\|_{L^\infty} \|v\|_{H^m} \|b \|_{H^m}+\sum_{|\a|\leq m}((b\cdot \nabla )D^\a v, D^\a
  b )_{L^2}.\n \\
  \eq
We observe that
  $$
 \sum_{|\a|\leq m}((b\cdot \nabla )D^\a b, D^\a v
   )_{L^2}=-\sum_{|\a|\leq m}((b\cdot \nabla )D^\a v, D^\a
  b )_{L^2},
  $$
 which is obvious by the integration by part. Thus, adding (\ref{lemm2}) to (\ref{lemm3}), we obtain
 \bb\label{lemm4}
 \frac12\frac{d}{dt} (\|v\|_{H^m} ^2+ \|b\|_{H^m} ^2 ) +\nu\|\nabla v\|_{H^m}^2
 \leq C(\|\nabla v\|_{L^\infty}  +\|\nabla
  b\|_{L^\infty} )(\|v\|_{H^m} ^2+ \|b\|_{H^m} ^2 ),
  \ee
  where we used the inequality, $ab\leq \frac12 (a^2 +b^2 )$.
 From (\ref{lemm4}) we first derive the
inequality
 \bq\label{lemm5}
\lefteqn{\|v(t)\|_{H^m} ^2+ \|b(t)\|_{H^m} ^2  +\nu\int_{t_0} ^t
\|\nabla v(s)\|_{H^m}^2ds}\hspace{.0in}\n \\
&& \leq (\|v(t_0)\|_{H^m} ^2+ \|b(t_0)\|_{H^m} ^2)\exp\left[
\int_{t_0} ^t( \|\nabla v(s)\|_{L^\infty}  +\|\nabla
  b(s)\|_{L^\infty})ds\right]\n \\
  \eq
 for all $0\leq t_0 <t $, which implies the continuation principle that
if
 $$ \int_{t_0} ^T( \|\nabla v(s)\|_{L^\infty}  +\|\nabla
  b(s)\|_{L^\infty})ds <\infty, $$
    then $\|v(T)\|_{H^m} + \|b(T)\|_{H^m}  <\infty$, and
we can continue our classical
  solution $(v(t), b(t))\in [H^m (\Bbb R^n)]^2$ up to $[t_0, T+\delta]$
  so that $(v,b)\in [C([0, T+\delta ]; H^m (\Bbb R^n ))]^2$ for some
  $\delta>0$.
Next, using  the estimate (\ref{lemm4}), we derive
 \bq
  \lefteqn{ \frac{d}{dt} \left\{(T-t)(\|v\|_{H^m} ^2+
\|b\|_{H^m} ^2 )\right\} +\nu (T-t)\|\nabla v\|_{H^m}^2
+(\|v\|_{H^m} ^2+ \|b\|_{H^m} ^2
)}\hspace{.5in}\n \\
&&\leq C_0(T-t)(\|\nabla v\|_{L^\infty}  +\|\nabla
  b\|_{L^\infty} )(\|v\|_{H^m} ^2+ \|b\|_{H^m} ^2 )
  \eq
  for a constant $C_0 =C_0 (m,n)$. We suppose
  $$
  \sup_{0<t<T} \left\{(T-t)(\|\nabla v(t)\|_{L^\infty}  +\|\nabla
  b(t)\|_{L^\infty} )\right\}<\frac{1}{2C_0}.
  $$
  Then,
  $$
\frac{d}{dt} \left\{(T-t)(\|v\|_{H^m} ^2+ \|b\|_{H^m} ^2 )\right\}
+\nu(T-t)\|\nabla v\|_{H^m}^2 +\frac12(\|v\|_{H^m} ^2+ \|b\|_{H^m}
^2) \leq 0,
$$
and integrating this over $[t_0, T]$, we have
 \bq
 \lefteqn{ \sup_{t_0 <t<T}
(T-t)(\|v\|_{H^m} ^2+ \|b\|_{H^m} ^2 )
  +\nu\int_{t_0}^T (T-t)\|\nabla v(t)\|_{H^m}^2dt}\hspace{.0in}\n \\
  &&+\frac12\int_{t_0}^T
  (\|v(t)\|_{H^m} ^2+ \|b(t)\|_{H^m} ^2)dt
  \leq (T-t_0)(\|v(t_0)\|_{H^m} ^2+ \|b(t_0)\|_{H^m} ^2.\n \\
 \eq
Since $H^m (\Bbb R^n) \hookrightarrow Lip (\Bbb R^n )$ for $m>n/2
+1$, we have
  \bqn
   \lefteqn{\int_{t_0}^T(\|\nabla v(t)\|_{L^\infty} + \|\nabla
b(t)\|_{L^\infty})dt \leq C\int_{t_0}^T(\|v(t)\|_{H^m}  + \|
b(t)\|_{H^m} )dt}\hspace{.0in}\\
&&\leq C\sqrt{T-t_0}\left[\int_{t_0}^T(\|v(t)\|_{H^m} ^2+
\|b(t)\|_{H^m} ^2)dt \right]^{\frac12}<\infty.
 \eqn
 Applying the continuation
principle derived above, we can
continue our local solution as described in the theorem. $\square$\\
\ \\
\noindent{\bf Proof of Theorem 1.2} We just observe that
$$
 (T-t)\|\nabla v (t)\|_{L^\infty}=
\|{V}\|_{L^\infty},\quad
  (T-t)\|\nabla b(t)\|_{L^\infty}=
\|\nabla {B}\|_{L^\infty}
 $$
for all $t\in (t_0, T)$. Hence, our smallness condition, $
 \|\nabla{V}\|_{L^\infty}
 +\|\nabla {B}\|_{L^\infty}<\eta,
 $
leads to
$$\sup_{t_0<t<T}(T-t)\left\{ \|\nabla v(t)\|_{{L^\infty}}
 +\|\nabla b(t)\|_{{L^\infty}}\right\}
<\eta.
 $$
 Applying Lemma 2.1, for initial time  at $t=t_0$, we
conclude  that $(v, b)\in [C ([t_0, T); H^m (\Bbb R^n ))]^2$
cannot have singularity at $t=t_0$, hence we need to have $V=B=0$. $\square$\\
\ \\
\noindent{\bf Proof of Theorem 1.3} We change variables from the
physical ones $(x,t) \in \Bbb R^n \times [0,T)$ to the `self-similar
variables' $(y,s)\in \Bbb R^n\times [0, \infty)$ as follows:
$$
y=\frac{x}{\sqrt{T-t}}, \quad s=\frac{1}{2} \log \left(
\frac{T}{T-t}\right).
$$
Based on this change of variables, we transform the functions
$(v,p)\mapsto (V, P)$ according to
 \bq
 v(x,t)&=&\frac{1}{\sqrt{T-t}} V(y,s ),\\
 b(x,t)&=&\frac{1}{\sqrt{T-t}} B(y,s ),\\
 p (x,t)&=&\frac{1}{\sqrt{T-t}}P (y,s ).
 \eq
 Substituting $(v,b, p )$ into (\ref{1.1})-(\ref{1.4}), we obtain the
 following equivalent evolution equations for
 $(V,P)$,
 \bb\label{selfmhd}
 \left\{ \aligned
  &  \frac{1 }{2}V_s +\frac{1}{2} V +\frac{1}{2}(y \cdot \nabla)V
  + (V\cdot \nabla )V =\nu \Delta V+(B\cdot \nabla )B-\nabla
( P +\frac12 |V|^2),\\
  &\frac{1 }{2}B_s +\frac{1}{2} B +\frac{1}{2}(y \cdot
  \nabla)B
  + (V\cdot \nabla )B =(B\cdot \nabla )V,\\
 & \qquad\mathrm{div}\, V=\mathrm{div}\, B=0,\\
 & V(y,0)=V_0 (y)=\sqrt{T} v_0 (\sqrt{T}y),\quad
  B(y,0)=B_0 (y)=\sqrt{T} b_0 (\sqrt{T}y).
  \endaligned \right.
 \ee
 In terms of $(V,B)$ the conditions
 (\ref{th12}) and
 (\ref{th122}) are translated
into
$$
 \sup_{0<s<\infty} \left(\|V(\cdot ,s)-\bar{V}(\cdot)\|_{L^1}
 +\|B(\cdot ,s)-\bar{B}(\cdot)\|_{L^1}\right) <\infty,
$$
and
 $$
 \lim_{s\to \infty}\|\nabla V(\cdot ,s) -\nabla\bar{V}(\cdot )\|_{L^\infty}=
 \lim_{s\to \infty}\|\nabla B(\cdot ,s) -\nabla \bar{B}(\cdot
 )\|_{L^\infty}=0,
 $$
 respectively, from which, thanks to the standard interpolation,
 we can have
 \bb\label{econv2}
  \lim_{s\to \infty}\| V(\cdot ,s) -\bar{V}(\cdot
   )\|_{H^1(B_R)}=
   \lim_{s\to \infty}\| B(\cdot ,s) -\bar{B}(\cdot
   )\|_{H^1(B_R)}=0
  \ee
  for all $0<R<\infty$, where  $B_R=\{ x\in \Bbb R^n \, |\, |x|<R\}$.
 Similarly to \cite{hou}, we consider scalar test functions
$\xi \in C^1_c (0,1)$ with $\int_0 ^1\xi (s)ds\neq 0$, $\psi \in
C^1_c (\Bbb R^n )$ and the vector test function $\phi =(\phi_1 ,
\cdots, \phi_n )\in C_c^1 (\Bbb
 R^n )$ with div $\phi=0$. We multiply the first equation of
 (\ref{selfmhd}) by $\xi
 (s-k)\phi (y)$, and integrate it over $\Bbb R^n\times [k, k+1]$,
 and then we integrate by part for the terms including the
 time derivative and the pressure term to obtain
 \bq\label{keq1}
 &&-\int_0^{1}\int_{\Bbb R^n} \xi _s(s) \phi(y)\cdot V(y,s+k)
 dyds\n\\
 &&+\int_0 ^{1}\int_{\Bbb R^n}\xi (s)\phi(y) \cdot[ V +(y \cdot \nabla)V
 +2
 (V\cdot \nabla )V](y,s+k)  dyds\n\\
 &&\quad-2\int_0 ^{1}\int_{\Bbb R^n}\xi (s)\phi(y) \cdot(B\cdot \nabla )B](y,s+k)
 dyds\n \\
 &&\quad +2\nu \int_0 ^{1}\int_{\Bbb R^n}\xi (s)\nabla \phi(y)\cdot \nabla V(y,s+k)
 dyds=0,
 \eq
 and
\bq\label{keq2}
  &&-\int_0^{1}\int_{\Bbb R^n} \xi _s(s) \psi(y) B(y,s+k)
 dyds\n\\
 &&+\int_0 ^{1}\int_{\Bbb R^n}\xi (s)\psi(y) [ B +(y \cdot \nabla)B
 +2
 (V\cdot \nabla )B](y,s+k)  dyds\n \\
 &&\quad -2\int_0 ^{1}\int_{\Bbb R^n}\xi (s)\psi(y)(B\cdot \nabla )V](y,s+k)  dyds=0.
 \eq
Passing to the limit $k\to \infty$ in (\ref{keq1})-(\ref{keq2}),
using the
  convergence (\ref{econv2}), $\int_0 ^1\xi _s(s)ds=0$ and $\int_0 ^1\xi
  (s)ds\neq 0$,
  we  find that $\bar{V}, \bar{B}\in C^1(\Bbb R^n)$ satisfies
$$
\int_{\Bbb R^n} [ \bar{V} +(y \cdot \nabla)\bar{V} +2 (\bar{V}\cdot
\nabla
 )\bar{V} -2
(\bar{B}\cdot \nabla
 )\bar{B}]\cdot \phi dy+2\nu \int_{\Bbb R^n} \nabla  \bar{V}\cdot \nabla \phi dy=0,
 $$
 $$
\int_{\Bbb R^n} [\bar{B} +(y \cdot \nabla)\bar{B} +2 (\bar{V}\cdot
\nabla
 )\bar{B} -2
(\bar{B}\cdot \nabla
 )\bar{V}]\psi dy=0,
 $$
 for all vector test function $\phi \in C_c^1 (\Bbb
 R^n)$ with div $\phi=0$, and scalar test function $\psi \in C^1_c
(\Bbb R^n )$
 Hence, there exists a scalar function
 $\bar{P'}$, which can be written without loss of generality as $\bar{P'}=
 \bar{P}+ \frac12 |\bar{B}|^2$ for another scalar function $\bar{P}$,  such that
 \bq\label{eleray1}
  \bar{V}+ (y\cdot\nabla )\bar{ V}
 +2(\bar{V}\cdot \nabla )\bar{V} =2\nu \Delta
 \bar{V}+
 2(\bar{B}\cdot \nabla )\bar{B}-2\nabla (\bar{P}+ \frac12 |\bar{B}|^2),\n \\
 \eq
 and
 \bb\label{eleray2}
   \bar{B}+ (y\cdot\nabla )\bar{ B}
 +2(\bar{V}\cdot \nabla )\bar{B} =2(\bar{B}\cdot \nabla
 )\bar{V}.
 \ee
 On the other hand, we can pass $s\to \infty$ directly in the
incompressibility equations for $V$ and $B$ in (\ref{selfmhd}) to
have
 \bb\label{eleray3}
\mathrm{ div} \, \bar{V}=\mathrm{div}\, \bar{B}=0.
 \ee
The equations (\ref{eleray1})-(\ref{eleray3})  show that $(\bar{V},
\bar{B} )$ is a classical solution of (\ref{le1})-(\ref{le3}).
Since, by hypothesis, $(\bar{V},\bar{B})$ satisfies the condition
(i)-(iii) of Theorem 1.1, we can deduce $\bar{V}=\bar{B}=0$ by that
theorem.
 Hence,  (\ref{econv1}) leads to
 $$\lim_{s\to \infty} \|\nabla V
 (s)\|_{L^\infty}=\lim_{s\to \infty} \|\nabla B
 (s)\|_{L^\infty}=0.
 $$
 Thus, for $\eta >0$ given in Lemma 2.1,
 there exists $s_1>0$ such that
 $$\|
\nabla V(s_1 )\|_{L^\infty}
 +\|\nabla B  (s_1 )\|_{L^\infty}<\eta .
 $$
Let us set $t_1=T[1-e^{2s_1} ]$. Going back to the original physical
variables, we have
 $$
  (T-t_1)\|\nabla v (t_1)\|_{L^\infty} +
  (T-t_1)\|\nabla b (t_1)\|_{L^\infty}<\eta .
 $$
 Applying Lemma 2.1, we conclude the proof. $\square$

\end{document}